\documentclass{article}
\usepackage{amsmath}

\newtheorem{theorem}{Theorem}

\newtheorem{example}[theorem]{Example}

\input{tcilatex}

\begin{document}

\title{Finding an \textsf{ARMA(p,q)} model given its spectral density or its 
$\rho _{k}$}
\author{Jan Vrbik \\
Dept. of Math, Brock University, St. Catharines, Ontario, Canada}

\begin{abstract}
An ARMA model can be fully determined based on either its spectral density,
or its correlogram, i.e. a formula for computing the corresponding $k^{th}$
serial correlation for any integer $k.$ In this article we describe how to
find, given one of these three ways of specifying the model, the other two.
\end{abstract}

\date{}
\maketitle

An \textsf{ARMA(p,q)} model is defined by its $\alpha $ and $\gamma $
coefficients (and their numbers, $p$ and $q$ respectively) and the value of $%
\sigma $ of the $\mathcal{N}(0,\sigma )$ distribution of the independent $%
\varepsilon $'s, thus:%
\begin{equation}
X_{n}=\alpha _{1}X_{n-1}+...+\alpha _{p}X_{n-p}+\varepsilon _{n}+\gamma
_{1}\varepsilon _{n-1}+...+\gamma _{q}\varepsilon _{n-q}  \label{mod}
\end{equation}%
where the $\gamma $'s can be arbitrary, but the $\alpha $'s must be such
that the $\lambda $ roots of 
\begin{equation}
\lambda ^{p}=\alpha _{1}\lambda ^{p-1}+\alpha _{2}\lambda ^{p-2}+...+\alpha
_{p}  \label{char}
\end{equation}%
are all, in absolute value, smaller than $1$ (they may be complex, in which
case they must go in complex-conjugate pairs). We call this the $\alpha $-$%
\gamma $ \emph{representation} (not always unique - there may be more than
one equivalent solution in terms of the $\gamma $ coefficients and,
consequently, of $\sigma $). We also require (without a loss of generality)
that 
\begin{equation*}
\lambda ^{q}+\gamma _{1}\lambda ^{q-1}+\gamma _{2}\lambda ^{q-2}+...+\gamma
_{q}=0
\end{equation*}%
and (\ref{char}) don't have any common roots.

Note that (\ref{mod}) applies the \textsf{MA} `filter' to the $\varepsilon $
\thinspace sequence, which then replaces the original $\varepsilon _{n}$ of
the \textsf{AR} part of the model. We should mention that doing the reverse,
i.e. using the usual \textsf{AR} model to generate a sequence of say $Y_{n}$
values, and then applying the \textsf{MA} filter to this sequence results in
values of $X_{n}$ with the same asymptotic properties (i.e. the same
correlogram and spectrum) as (\ref{mod}). Note that throughout this article
we assume that the corresponding stochastic process has equilibrated, and
has therefore reached its stationary mode.\bigskip

The same model can also be identified by the formula for $\rho _{k},$ which
must have the following form:%
\begin{equation}
\rho _{k}=\sum_{i}P_{i}(k)\cdot \theta _{i}^{k}+\sum_{j}\left( Q_{j}(k)\cdot
\sin (a_{j}\cdot k)\cdot \theta _{j}^{k}\overset{}{+}R_{j}(k)\cdot \cos
(a_{j}\cdot k)\cdot \theta _{j}^{k}\right) \ \ \ \ \ \ \ \ \ \ \ \ \ \text{%
when }k>q-p  \label{corr}
\end{equation}%
where each $\theta $ (they are all \emph{real}) must be, in absolute value,
less than $1,$ $P_{i}(k)$, $Q_{j}(k)$ and $R_{j}(k)$ are \emph{polynomials}
in $k,$ and both summations are \emph{finite}. In addition, one also needs
to be given the value of the `exceptional' $\rho _{k}$s (for $k$ from $0$ to 
$q-p$), if any, and the value of $\sigma $ (or, equivalently, of the common $%
X_{n}$ variance $V$). This is called the \emph{correlogram} representation.
Note that the absolute-value condition is necessary but not sufficient to
guarantee that the formula corresponds to a legitimate and stationary 
\textsf{ARMA} model (the legitimacy is easily established in either the $%
\alpha $-$\gamma $ or the spectral representation).

Also note that (\ref{corr}) can be always simplified to 
\begin{equation}
\rho _{k}=\sum_{i}\tilde{P}_{i}(k)\cdot \tilde{\theta}_{i}^{k}  \label{corr2}
\end{equation}%
where some of the $\tilde{\theta}$ (and the corresponding $\tilde{P}$) may
now be complex (they must appear in complex-conjugate pairs). When expanded
(each power of $k$ becomes a separate term), the total number of terms in (%
\ref{corr2}) must equal to $p$.

It is easy to go from the `real' representation (\ref{corr}) to the
`complex' representation (\ref{corr2}), and reverse. In this article, we
prefer to use (\ref{corr2}).\bigskip

Finally, another way of fully specifying the model is to know its \emph{%
spectral density} 
\begin{equation}
\omega (\beta )=\frac{S(\cos \beta )}{T(\cos \beta )}  \label{spe}
\end{equation}%
which must be a rational function of $\cos \beta $ (both $S$ and $T$ are
polynomials, of degree $q$ and $p$ respectively, not sharing any roots),
together with the value of either $\sigma $ or $V$ - this is the \emph{%
spectral} representation. Note that $\int_{0}^{\pi }\omega (\beta )d\beta $
must equal to $\pi $ (making the \emph{average} value of the spectral
density equal to $1$) and $\omega (\beta )$ must be \emph{non-negative} and
bounded (no singularities) - this can be easily verified by plotting it. The
two conditions are necessary and \emph{sufficient} to make $\omega (\beta )$
`legitimate'.\bigskip

Now, let us deal with the following obvious issue: given one of the three
possible representations, how do we find the other two?

\paragraph{$\protect\alpha $-$\protect\gamma $ representation converted to
correlogram:}

We must first find the $X$ variance (the common variance of all the $X_{i}$%
s), and the $k^{th}$ serial correlation coefficient (let's call them $\tilde{%
V}$ and $\tilde{\rho}_{k}$ respectively) of the \textsf{AR} part of the
model by solving%
\begin{equation}
\tilde{\rho}_{i}=\sum_{j=1}^{p}\alpha _{j}\tilde{\rho}_{|i-j|}\ \ \ \ \ \ \
\ \ \ \ \ \ \ \ \ \ 1\leq i\leq p-1  \label{A1}
\end{equation}%
for $\tilde{\rho}_{1},\tilde{\rho}_{2},...\tilde{\rho}_{p-1}$ ($\tilde{\rho}%
_{0}$ is always equal to $1$), then setting%
\begin{equation}
\tilde{\rho}_{k}=\dsum_{i}\tsum_{j=1}^{m_{i}}A_{i,j}\cdot k^{j-1}\cdot
\lambda _{i}^{k}  \label{A2}
\end{equation}%
where $\lambda _{i}$ are the roots of (\ref{char}) and $m_{i}$ is the
multiplicity of $\lambda _{i}$ (note that the double summation must account
for exactly $p$ roots) and solving for the $A_{i,j}$ coefficients based on
the knowledge of the first $p$ values of $\tilde{\rho}_{k}$ (including $%
\tilde{\rho}_{0}=1$). Note that the resulting formula is correct for any $%
k>-p$. Also note that it is not always possible to find $\lambda _{i}$ in an
exact form (as we usually like to do), and we have to switch to decimals
(posing extra computational challenges, such as the difficulty of
establishing whether two or more roots are identical or not). In this
article, we will bypass this problem by working with exact solutions only
(justified when $p<5$).

Having these, we find%
\begin{equation}
\tilde{V}=\frac{\sigma ^{2}}{1-\sum_{j=1}^{p}\alpha _{j}\tilde{\rho}_{j}}
\label{A3}
\end{equation}%
(the common variance of the $Y_{n}$ sequence).

We can now proceed to compute the variance and serial correlation of the
complete model by:%
\begin{eqnarray}
V &=&\hat{V}\cdot \sum_{i,j=0}^{q}\gamma _{i}\gamma _{j}\hat{\rho}_{|i-j|}
\label{V} \\
\rho _{k} &=&\frac{\sum_{i,j=0}^{q}\gamma _{i}\gamma _{j}\hat{\rho}_{k+i-j}}{%
\sum_{i,j=0}^{q}\gamma _{i}\gamma _{j}\hat{\rho}_{|i-j|}}\ \ \ \ \ \ \ \ 
\text{when }k>q-p  \label{V1} \\
\rho _{k} &=&\frac{\sum_{i,j=0}^{q}\gamma _{i}\gamma _{j}\hat{\rho}_{|k+i-j|}%
}{\sum_{i,j=0}^{q}\gamma _{i}\gamma _{j}\hat{\rho}_{|i-j|}}\ \ \ \ \ \ \text{%
for the remaining non-negative }k\ \text{(if any)}  \label{V2}
\end{eqnarray}%
with the understanding that $\gamma _{0}=1.$

To simplify the resulting expressions, we can use Maple's `expand' and
`convert($\cdots $,power)' commands.

\paragraph{Correlogram to spectral:}

This can be achieved rather easily by%
\begin{equation*}
\omega (\beta )=1+2\sum_{k=1}^{\infty }\rho _{k}\cdot \cos (k\beta )
\end{equation*}%
and simplifying the answer. For this step to work well, it is essential to
use an \emph{exact} (no decimals) formula for $\rho _{k}.$

\paragraph{Spectral to $\protect\alpha $-$\protect\gamma $ :}

First we have to find, \emph{separately}, the roots of the $S$ and $T$
polynomials from (\ref{spe}), and for each of these (let's call them,
generically, $\theta $), to solve the following quadratic equation for $y$ 
\begin{equation}
y+\frac{1}{y}=2\theta  \label{the}
\end{equation}%
getting two reciprocal solutions. For the \textsf{AR} part (represented by
the denominator $T$) one of these must always be (in absolute value) smaller
than $1$ (the other is bigger than $1$). Collecting the \textsf{smaller}
solutions (let's call them $y_{1},$ $y_{2},$ ...$y_{p}$), we expand%
\begin{equation}
(y-y_{1})(y-y_{2})...(y-y_{p})  \label{exp}
\end{equation}%
The resulting coefficients of $y^{p-1},$ $y^{p-2},....$ $y^{1}$ and $y^{0}$
yield $-\alpha _{1},$ $-\alpha _{2},..-\alpha _{p-1}$ and $-\alpha _{p}$
(respectively) of the \textsf{AR} part of the model.

To get the $\gamma $ coefficients (of the \textsf{MA} part), we do the same
thing with the roots of the numerator $S$, except now

\begin{itemize}
\item some of the $y$'s may have the absolute value of $1,$

\item we can select \emph{either one} of the two reciprocal solutions (it
does not have to be the `smaller' one),

\item selecting a \emph{complex} $y,$ we must also select its complex \emph{%
conjugate} (from yet \emph{another} pair of reciprocal solutions),

\item after expanding $(y-y_{1})(y-y_{2})...(y-y_{q})$, the coefficients
yield (in the same manner as before) $+\gamma _{1},$ $+\gamma _{2},..+\gamma
_{p-1}$ and $+\gamma _{q}$.\bigskip
\end{itemize}

Note that we have just completed a full cycle through the $\alpha $-$\gamma $
$\longrightarrow $ correlogram $\longrightarrow $ spectral $\longrightarrow $
$\alpha $-$\gamma $ `triangle', which would be sufficient to get from one
representation to any of the other two. For completeness, we continue
investigating the remaining three direct `links', namely:

\paragraph{$\protect\alpha $-$\protect\gamma $ to spectral:}

This is again done quite easily by simplifying%
\begin{equation}
\omega (\beta )=\frac{\sigma ^{2}}{V}\cdot \frac{\sum_{j=0}^{q}\gamma
_{j}e^{i\beta j}\cdot \sum_{j=0}^{q}\gamma _{j}e^{-i\beta j}}{%
\sum_{j=0}^{p}\alpha _{j}e^{i\beta j}\cdot \sum_{j=0}^{p}\alpha
_{j}e^{-i\beta j}}  \label{ome}
\end{equation}%
(with the understanding that $\gamma _{0}=1$ and $\alpha _{0}=-1$), i.e.
converting it to a function of $\cos (\beta )$ by Maple$^{\prime }$s
`convert($\cdots $,trig)'. Note that $\int_{0}^{\pi }\omega (\beta )~d\beta $
must equal to $\pi $, which simplifies finding the $\frac{\sigma ^{2}}{V}$
factor. Also not that the previous formula implies that%
\begin{equation}
\omega (0)=\frac{\sigma ^{2}}{V}\cdot \left( \frac{\sum_{j=0}^{q}\gamma _{j}%
}{\sum_{j=0}^{p}\alpha _{j}}\right) ^{2}  \label{zer}
\end{equation}

\paragraph{Correlogram to $\protect\alpha $-$\protect\gamma $ :}

Based on the $\rho _{k}$ formula (valid for $k>m,$ where $m$ could be
negative - verify its exact value), we can easily see what the roots of (\ref%
{char}) are, including their multiplicity (based on the degree of the $%
\tilde{P}_{i}$ polynomials). Expanding%
\begin{equation*}
(\lambda -\lambda _{1})(\lambda -\lambda _{2})...(\lambda -\lambda _{p})
\end{equation*}%
yields the $\alpha $'s, and correspondingly the $\tilde{\rho}_{k}$ formula
(as explained in the `$\alpha $-$\gamma $ to correlogram' section). Since $%
m=q-p$ (where $p$ is known already) we can find $q.$ We then solve the $q-p$
`special' equations (\ref{V2}) and the $p$ equations we get from each term
of (\ref{V1}) for $\gamma _{1},$ $\gamma _{2},...\gamma _{q}$ (when $p>q,$
these $p$ equations for $q$ unknowns are guaranteed to have at least one
solution, but one must work with \emph{exact} quantities only - no
decimals). Note that `terms' of (\ref{V1}) are defined in the sense of \emph{%
expanded} (\ref{corr2}). Also note that the resulting equations are of
non-linear (multinomial) type - see our Examples; unless a special technique
is developed, solving them beyond the case of $q=3$ becomes practically
impossible (one can then use the `correlogram $\longrightarrow $ spectral $%
\longrightarrow $ $\alpha $-$\gamma $' route, which works fine).

\paragraph{Spectral to correlogram:}

The basic formula (from the Fouries-series theory) is%
\begin{equation}
\rho _{k}=\frac{1}{2\pi }\int_{0}^{2\pi }\omega (\beta )\cos (k\beta )d\beta
\label{fur}
\end{equation}%
Unfortunately, Maple cannot deal with a general integer $k$ (yet?), but we
can do this on our own, using contour integration. This results in the
following algorithm (see the appendix):

\begin{itemize}
\item Replace, in (\ref{spe}), $\cos \beta $ by%
\begin{equation}
\frac{y+\frac{1}{y}}{2}  \label{cos}
\end{equation}%
and find the corresponding (\emph{complex}) partial-fraction expansion (in $%
y $). The result will be a finite sum of terms of the following two types,
either%
\begin{equation*}
\frac{c}{(y-\lambda )^{\ell }}
\end{equation*}%
where $\ell $ is a positive integer (usually equal to $1$) and $c$ and $%
\lambda $ are specific (potentially complex) constants, or%
\begin{equation*}
c\cdot y^{m}
\end{equation*}%
where $m$ can be any (positive, negative, or zero) integer.

\item Ignore (delete) terms of the former type with $|\lambda |>1,$ and
terms of the latter type with $m$ negative.

\item Replace each remaining term by 
\begin{equation*}
c\cdot \binom{k-1}{k-\ell }\lambda ^{k-\ell }
\end{equation*}%
in the first case (the sum of these constitutes the general formula), and by%
\begin{equation*}
c\cdot \delta _{m,k}
\end{equation*}%
in the second case (these contribute to the `special' cases of $\rho _{k}$%
.), where $\delta _{m,k}$ is the Kronecker's delta.
\end{itemize}

\begin{example}
Let us consider the following \textsf{ARMA(3,2)} model:%
\begin{equation*}
X_{n}=\frac{133}{60}X_{n-1}-\frac{49}{30}X_{n-2}+\frac{2}{5}%
X_{n-3}+\varepsilon _{n}-4\varepsilon _{n-1}+5\varepsilon _{n-2}
\end{equation*}%
with $\sigma =\frac{1}{10}.$
\end{example}

\subparagraph{$\protect\alpha $-$\protect\gamma $ $\longrightarrow $
correlogram:}

Using (\ref{A1}), (\ref{A2}) and (\ref{A3}), we find%
\begin{eqnarray}
\hat{\rho}_{k} &=&\frac{125\cdot \left( \frac{4}{5}\right) ^{k}-135\cdot
\left( \frac{3}{4}\right) ^{k}+27\cdot \left( \frac{2}{3}\right) ^{k}}{17}
\label{arp} \\
\hat{V} &=&\frac{68}{35}  \notag
\end{eqnarray}%
(note that the former formula is good for all $k>-3$), and%
\begin{eqnarray}
\rho _{k} &=&\frac{1525\cdot \left( \frac{4}{5}\right) ^{k}-1599\cdot \left( 
\frac{3}{4}\right) ^{k}+300\cdot \left( \frac{2}{3}\right) ^{k}}{226}
\label{rho} \\
V &=&\frac{113}{14}  \notag
\end{eqnarray}%
(good for all non-negative $k$).

\subparagraph{$\protect\alpha $-$\protect\gamma $ $\longrightarrow $
spectral:}

Using (\ref{ome}), we get%
\begin{eqnarray}
\omega (\beta ) &=&\frac{\frac{\sigma ^{2}}{V}\cdot \left( 1-4e^{i\beta
}+5e^{2i\beta }\right) \cdot \left( 1-4e^{-i\beta }+5e^{-2i\beta }\right) }{%
\left( -1+\frac{133}{60}e^{i\beta }-\frac{49}{30}e^{2i\beta }+\frac{2}{5}%
e^{3i\beta }\right) \cdot \left( -1+\frac{133}{60}e^{-i\beta }-\frac{49}{30}%
e^{-2i\beta }+\frac{2}{5}e^{-3i\beta }\right) }  \notag \\
&=&\frac{2016}{113}\cdot \frac{8-12\cos \beta +5\cos ^{2}\beta }{%
13325-38092\cos \beta +36288\cos ^{2}\beta -11520\cos ^{3}\beta }
\label{sec}
\end{eqnarray}%
after proper simplification. Note that 
\begin{equation*}
\frac{\sigma ^{2}}{V}=\frac{7}{5650}
\end{equation*}%
based on the $\int_{0}^{\pi }\omega (\beta )~d\beta =\pi $ condition.

\subparagraph{Spectral $\longrightarrow ~\protect\alpha $-$\protect\gamma $:}

The roots of%
\begin{equation*}
8-12\theta +5\theta ^{2}=0
\end{equation*}%
are $\theta =\dfrac{6+2i}{5}$ and $\dfrac{6-2i}{5}.$ Solving (\ref{the}) for 
$y$ yields: $2+i$ and $\dfrac{2-i}{5}$ for the first $\theta ,$ and $2-i$
and $\dfrac{2+i}{5}$ for the second one. Selecting one value from each pair
(note that they must go in complex-conjugates), (\ref{exp}) yields%
\begin{equation*}
\left( y-2-i\right) \left( y-2+i\right) =y^{2}-4y+5
\end{equation*}%
verifying the \textsf{MA(2)} part of the model. Note that the second
possible choice would yield%
\begin{equation*}
\left( y-\frac{2-i}{5}\right) \left( y-\frac{2+i}{5}\right) =y^{2}-\frac{4}{5%
}y+\frac{1}{5}
\end{equation*}%
resulting in a model not identical but \emph{equivalent} to the original one
(after adjusting $\sigma $ accordingly, this time making it equal to $\frac{1%
}{2}$).

The roots of the denominator are $\theta =\dfrac{13}{12},$ $\dfrac{25}{24}$
and $\dfrac{41}{40}.$ Solving for $y$ yields: $\dfrac{3}{2}$ and $\dfrac{2}{3%
}$ for the first $\theta ,$ $\dfrac{4}{3}$ and $\dfrac{3}{4}$ for the second
one, and $\dfrac{5}{4}$ and $\dfrac{4}{5}$ for the last one Taking the
smaller value of each pair, we get%
\begin{equation}
\left( y-\frac{2}{3}\right) \left( y-\frac{3}{4}\right) \left( y-\frac{4}{5}%
\right) =y^{3}-\frac{133}{60}y^{2}+\frac{49}{30}y-\frac{2}{5}  \label{alp}
\end{equation}%
thus recovering the \textsf{AR(3)} part of the model.

We get $V,$ based on (\ref{zer}), by%
\begin{equation*}
V=\frac{\sigma ^{2}\cdot (1-4+5)^{2}}{\omega (0)\cdot \left( 1-\frac{133}{60}%
+\frac{49}{30}-\frac{2}{5}\right) ^{2}}=\frac{113}{14}
\end{equation*}%
(check).

\subparagraph{Spectral$~\longrightarrow ~$correlogram:}

Replacing $\cos \beta $ by (\ref{cos}), $\omega (\beta )$ can be expressed
in a partial-fraction form as%
\begin{equation*}
-\frac{7625}{904\left( y-\frac{5}{4}\right) }-\frac{225}{113\left( y-\frac{3%
}{2}\right) }+\frac{1066}{113\left( y-\frac{4}{3}\right) }+\frac{610}{%
113\left( y-\frac{4}{5}\right) }+\frac{100}{113\left( y-\frac{2}{3}\right) }-%
\frac{4797}{904\left( y-\frac{3}{4}\right) }
\end{equation*}%
Based on the last three terms (ignoring the first three), we get%
\begin{eqnarray}
\rho _{k} &=&\frac{610}{113}\cdot \left( \frac{4}{5}\right) ^{k-1}-\frac{4797%
}{904}\cdot \left( \frac{3}{4}\right) ^{k-1}+\frac{100}{113}\cdot \left( 
\frac{2}{3}\right) ^{k-1}  \notag \\
&=&\frac{1525}{226}\cdot \left( \frac{4}{5}\right) ^{k}-\frac{1599}{226}%
\cdot \left( \frac{3}{4}\right) ^{k}+\frac{300}{226}\cdot \left( \frac{2}{3}%
\right) ^{k}  \label{Rk}
\end{eqnarray}%
good for all non-negative $k$. This agrees with (\ref{rho}).

\subparagraph{Correlogram$~\longrightarrow ~$spectral$:$}

Simplifying%
\begin{equation*}
1+\frac{1525}{113}\sum_{k=1}^{\infty }\left( \frac{4}{5}\right) ^{k}\cos
(k\cdot \beta )-\frac{1599}{113}\sum_{k=1}^{\infty }\left( \frac{3}{4}%
\right) ^{k}\cos (k\cdot \beta )+\frac{300}{113}\sum_{k=1}^{\infty }\left( 
\frac{2}{3}\right) ^{k}\cos (k\cdot \beta )
\end{equation*}%
yields (\ref{sec}).

\subparagraph{Correlogram $\longrightarrow $ $\protect\alpha $-$\protect%
\gamma $ :}

From (\ref{rho}) it is immediately obvious that there are three simple $%
\lambda $ roots, namely $\frac{4}{5},$ $\frac{2}{3}$ and $\frac{3}{4}.$ (\ref%
{alp}) then yields the value of the $\alpha $ parameters. Knowing these, we
can find (\ref{arp}).

Since $p=3$ and, based on the $k$ range of (\ref{Rk}), $q-p=-1,$ $q$ must be
equal to $2$.

We now need to solve the following three nonlinear equations for $\gamma
_{1} $ and $\gamma _{2}$, based on (\ref{V1}) and (\ref{Rk}):%
\begin{eqnarray*}
&&\frac{125}{17}\cdot \sum_{i,j=0}^{2}\left( \frac{4}{5}\right) ^{i-j}\cdot
\gamma _{i}\gamma _{j}\overset{}{=}\frac{1525}{226}\cdot \sum_{i,j=0}^{2}%
\hat{\rho}_{|i-j|}\cdot \gamma _{i}\gamma _{j} \\
&&-\frac{135}{17}\cdot \sum_{i,j=0}^{2}\left( \frac{3}{4}\right) ^{i-j}\cdot
\gamma _{i}\gamma _{j}\overset{}{=}-\frac{1599}{226}\cdot \sum_{i,j=0}^{2}%
\hat{\rho}_{|i-j|}\cdot \gamma _{i}\gamma _{j} \\
&&\frac{27}{17}\cdot \sum_{i,j=0}^{2}\left( \frac{2}{3}\right) ^{i-j}\cdot
\gamma _{i}\gamma _{j}\overset{}{=}\frac{300}{226}\cdot \sum_{i,j=0}^{2}\hat{%
\rho}_{|i-j|}\cdot \gamma _{i}\gamma _{j}
\end{eqnarray*}%
where $\hat{\rho}_{k}$ is defined in (\ref{arp}), and with the understanding
that $\gamma _{0}=1$. These can be reduced to%
\begin{eqnarray*}
334(1+\gamma _{1}^{2}+\gamma _{2}^{2})+1057\gamma _{1}(1+\gamma
_{2})+2268\gamma _{2} &=&0 \\
1109(1+\gamma _{1}^{2}+\gamma _{2}^{2})+3332\gamma _{1}(1+\gamma
_{2})+6678\gamma _{2} &=&0 \\
155(1+\gamma _{1}^{2}+\gamma _{2}^{2})+455\gamma _{1}(1+\gamma
_{2})+882\gamma _{2} &=&0
\end{eqnarray*}%
Even though seemingly over-determined, they must have at least one (two in
this case) real solution, namely $\gamma _{1}=-4$ and $\gamma _{2}=5$
(check), but also $\gamma _{1}=-\frac{4}{5}$ and $\gamma _{2}=\frac{1}{5}$
(note that one can get this solution from the previous one, namely $\mathbf{%
\gamma }=(1,-4,5)$, by reversing the order and dividing by the $5,$ the
originally last $\gamma $).

\begin{example}
Let us consider the following \textsf{ARMA(2,3)} model:%
\begin{equation*}
X_{n}=X_{n-1}-\frac{1}{2}X_{n-2}+\varepsilon _{n}+3\varepsilon
_{n-1}+3\varepsilon _{n-2}+\varepsilon _{n-3}
\end{equation*}%
with $\sigma =\frac{1}{10}.$
\end{example}

\subparagraph{$\protect\alpha $-$\protect\gamma $ $\longrightarrow $
correlogram:}

Using (\ref{A1}), (\ref{A2}) and (\ref{A3}), we find%
\begin{eqnarray}
\hat{\rho}_{k} &=&\left( \frac{1}{2}+\frac{i}{6}\right) \left( \frac{1-i}{2}%
\right) ^{k}+\left( \frac{1}{2}-\frac{i}{6}\right) \left( \frac{1+i}{2}%
\right) ^{k}  \label{arp2} \\
\hat{V} &=&\frac{3}{125}  \notag
\end{eqnarray}%
(note that the former formula is good for all $k>-2$), and%
\begin{eqnarray}
\rho _{k} &=&\frac{41+38i}{100}\cdot \left( \frac{1-i}{2}\right) ^{k}+\frac{%
41-38i}{100}\cdot \left( \frac{1+i}{2}\right) ^{k}\ \ \ \ \ \ \ \text{when }%
k>1  \label{rho2} \\
\rho _{0} &=&1\ \ \ \ \text{and \ \ \ }\rho _{1}=\frac{81}{100}  \notag \\
V &=&1  \notag
\end{eqnarray}%
Note that the general formula can be also written in an explicitly real form
of%
\begin{equation*}
\rho _{k}=\left( \frac{1}{\sqrt{2}}\right) ^{k}\cdot \left( \frac{38\sin 
\frac{k\ \pi }{4}+41\cos \frac{k\ \pi }{4}}{50}\right)
\end{equation*}

\subparagraph{$\protect\alpha $-$\protect\gamma $ $\longrightarrow $
spectral:}

Using (\ref{ome}), we get%
\begin{eqnarray}
\omega (\beta ) &=&\frac{\frac{\sigma ^{2}}{V}\cdot \left( 1+3e^{i\beta
}+3e^{2i\beta }+e^{3i\beta }\right) \cdot \left( 1+3e^{-i\beta
}+3e^{-2i\beta }+e^{-3i\beta }\right) }{\left( -1+e^{i\beta }-\frac{1}{2}%
e^{2i\beta }\right) \cdot \left( -1+e^{-i\beta }-\frac{1}{2}e^{-2i\beta
}\right) }  \notag \\
&=&\frac{8}{25}\cdot \frac{1+3\cos \beta +3\cos ^{2}\beta +\cos ^{3}\beta }{%
5-12\cos \beta +8\cos ^{2}\beta }  \label{sec2}
\end{eqnarray}%
after proper simplification. Note that 
\begin{equation*}
\frac{\sigma ^{2}}{V}=\frac{1}{100}
\end{equation*}%
based on the $\int_{0}^{\pi }\omega (\beta )~d\beta =\pi $ condition.

\subparagraph{Spectral $\longrightarrow ~\protect\alpha $-$\protect\gamma $:}

The three roots of%
\begin{equation*}
1+3\theta +3\theta ^{2}+\theta ^{3}=0
\end{equation*}%
are all equal to $-1$. Solving (\ref{the}) for $y$ always yields a double
root of $-1.$ Selecting one value from each pair (not much choice now), (\ref%
{exp}) expands to%
\begin{equation*}
\left( y+1\right) \left( y+1\right) (y+1)=y^{3}+3y^{2}+3y+1
\end{equation*}%
verifying the \textsf{MA(2)} part of the model.

The roots of the denominator are $\theta =\frac{3+i}{4},$ and $\frac{3-i}{4}%
. $ Solving for $y$ yields: $1+i$ and $\dfrac{1-i}{2}$ for the first $\theta 
$ and $1-i$ and $\dfrac{1+i}{2}$ for the second one Taking the smaller value
of each pair, we get%
\begin{equation}
\left( y-\frac{1-i}{2}\right) \left( y-\frac{1+i}{2}\right) =y^{2}-y+\frac{1%
}{2}  \label{alp2}
\end{equation}
thus recovering the \textsf{AR(3)} part of the model.

We get $V,$ based on (\ref{zer}), by%
\begin{equation*}
V=\frac{\sigma ^{2}\cdot (1+3+3+1)^{2}}{\omega (0)\cdot \left( 1-1+\frac{1}{2%
}\right) ^{2}}=1
\end{equation*}%
(check).

\subparagraph{Spectral$~\longrightarrow ~$correlogram:}

Replacing $\cos \beta $ by (\ref{cos}), $\omega (\beta )$ can be expressed
in a partial-fraction form as%
\begin{equation*}
\frac{9}{50}+\frac{y}{50}+\frac{1}{50y}+\frac{3+79i}{100\cdot (y-1-i)}+\frac{%
3-79i}{100\cdot (y-1+i)}+\frac{79+3i}{200\cdot \left( y-\frac{1+i}{2}\right) 
}+\frac{79-3i}{200\cdot \left( y-\frac{1-i}{2}\right) }
\end{equation*}%
Based on the first two and the last two terms (ignoring the rest), we get%
\begin{equation}
\rho _{k}=\frac{9}{50}\delta _{0k}+\frac{1}{50}\delta _{1k}+\frac{41+38i}{100%
}\cdot \left( \frac{1-i}{2}\right) ^{k}+\frac{41-38i}{100}\cdot \left( \frac{%
1+i}{2}\right) ^{k}  \label{Rk2}
\end{equation}%
which agrees with (\ref{rho2}).

\subparagraph{Correlogram$~\longrightarrow ~$spectral$:$}

Simplifying%
\begin{equation*}
1+\frac{81}{50}\cos \beta +\frac{41+38i}{50}\sum_{k=2}^{\infty }\left( \frac{%
1-i}{2}\right) ^{k}\cos (k\cdot \beta )+\frac{41-38i}{50}\sum_{k=2}^{\infty
}\left( \frac{1+i}{2}\right) ^{k}\cos (k\cdot \beta )
\end{equation*}
yields (\ref{sec2}).

\subparagraph{Correlogram $\longrightarrow $ $\protect\alpha $-$\protect%
\gamma $ :}

From (\ref{rho2}) it is immediately obvious that there are two simple $%
\lambda $ roots, namely $\frac{1-i}{2}$ and $\frac{1+i}{2}.$ (\ref{alp2})
then yields the value of the $\alpha $ parameters. Knowing these, we can
find (\ref{arp2}).

Since $p=2$ and, based on the $k$ range of the general part of (\ref{Rk2}), $%
q-p=1,$ $q$ must be equal to $3$.

We now need to solve the following three nonlinear equations for $\gamma
_{1},$ $\gamma _{2}$ and $\gamma _{3}$, based on (\ref{V2}), (\ref{V1}) and (%
\ref{Rk2}):%
\begin{eqnarray*}
&&\sum_{j,\ell =0}^{3}\hat{\rho}_{|1+j-\ell |}\cdot \gamma _{j}\gamma _{\ell
}\overset{}{=}\frac{81}{100}\cdot \sum_{j,\ell =0}^{3}\hat{\rho}_{|j-\ell
|}\cdot \gamma _{j}\gamma _{\ell } \\
&&\left( \frac{1}{2}+\frac{i}{6}\right) \cdot \sum_{j,\ell =0}^{3}\left( 
\frac{1-i}{2}\right) ^{j-\ell }\cdot \gamma _{i}\gamma _{\ell }\overset{}{=}%
\frac{41+38i}{100}\cdot \sum_{j,\ell =0}^{3}\hat{\rho}_{|j-\ell |}\cdot
\gamma _{j}\gamma _{\ell } \\
&&\left( \frac{1}{2}-\frac{i}{6}\right) \cdot \sum_{j,\ell =0}^{3}\left( 
\frac{1+i}{2}\right) ^{j-\ell }\cdot \gamma _{i}\gamma _{\ell }\overset{}{=}%
\frac{41-38i}{100}\cdot \sum_{j,\ell =0}^{3}\hat{\rho}_{|j-\ell |}\cdot
\gamma _{j}\gamma _{\ell }
\end{eqnarray*}%
where $\hat{\rho}_{k}$ is defined in (\ref{arp2}), and with the
understanding that $\gamma _{0}=1$. These can be reduced to%
\begin{eqnarray*}
&&43(1+\gamma _{1}^{2}+\gamma _{2}^{2}+\gamma _{3}^{2})-26(\gamma
_{1}+\gamma _{1}\gamma _{2}+\gamma _{2}\gamma _{3})-69(\gamma _{2}+\gamma
_{1}\gamma _{3})-56\gamma _{3}\overset{}{=}0 \\
&&(7+12i)(1+\gamma _{1}^{2}+\gamma _{2}^{2}+\gamma _{3}^{2})-(4-6i)(\gamma
_{1}+\gamma _{1}\gamma _{2}+\gamma _{2}\gamma _{3}) \\
&&\ \ \ \ \ \ \ \ \ \ \ \ \ \ \ \ \ \ \ \ \ \ \ \ \ \ \ \ \ \ \ \ \ \ \ \ \
\ \ \ -(16+41i)(\gamma _{2}+\gamma _{1}\gamma _{3})+(16-84i)\gamma _{3}%
\overset{}{=}0
\end{eqnarray*}%
and the corresponding complex conjugate. This time, there is only one
solution, namely $\gamma _{1}=3,$ $\gamma _{2}=3$ and $\gamma _{3}=1$
(check).

\begin{example}
Let us consider the following \textsf{ARMA(3,2)} model:%
\begin{equation*}
X_{n}=\frac{3}{2}X_{n-1}-\frac{3}{4}X_{n-2}+\frac{1}{8}X_{n-3}+\varepsilon
_{n}-2\varepsilon _{n-1}+2\varepsilon _{n-2}
\end{equation*}%
with $\sigma =\frac{1}{10}.$
\end{example}

\subparagraph{$\protect\alpha $-$\protect\gamma $ $\longrightarrow $
correlogram:}

Using (\ref{A1}), (\ref{A2}) and (\ref{A3}), we find%
\begin{eqnarray}
\hat{\rho}_{k} &=&\left( 1+\frac{15}{22}\cdot k+\frac{3}{22}\cdot
k^{2}\right) \left( \frac{1}{2}\right) ^{k}  \label{arp3} \\
\hat{V} &=&\frac{176}{2025}  \notag
\end{eqnarray}%
(note that the former formula is good for all $k>-3$), and%
\begin{eqnarray}
\rho _{k} &=&\left( 1+\frac{3}{44}\cdot k+\frac{15}{44}\cdot k^{2}\right)
\left( \frac{1}{2}\right) ^{k}\ \ \ \ \ \ \text{when }k>-1  \label{rhox} \\
V &=&\frac{176}{2025}  \notag
\end{eqnarray}

\subparagraph{$\protect\alpha $-$\protect\gamma $ $\longrightarrow $
spectral:}

Using (\ref{ome}), we get%
\begin{eqnarray}
\omega (\beta ) &=&\frac{\frac{\sigma ^{2}}{V}\cdot \left( 1-2e^{i\beta
}+2e^{2i\beta }\right) \cdot \left( 1-2e^{-i\beta }+2e^{-2i\beta }\right) }{%
\left( -1+\frac{3}{2}e^{i\beta }-\frac{3}{4}e^{2i\beta }+\frac{1}{8}%
e^{3i\beta }\right) \cdot \left( -1+\frac{3}{2}e^{-i\beta }-\frac{3}{4}%
e^{-2i\beta }+\frac{1}{8}e^{-3i\beta }\right) }  \notag \\
&=&\frac{81}{11}\cdot \frac{5-12\cos \beta +8\cos ^{2}\beta }{\left( 5-4\cos
\beta \right) ^{3}}  \label{sec3}
\end{eqnarray}%
after proper simplification. Note that 
\begin{equation*}
\frac{\sigma ^{2}}{V}=\frac{81}{704}
\end{equation*}%
based on the $\int_{0}^{\pi }\omega (\beta )~d\beta =\pi $ condition.

\subparagraph{Spectral $\longrightarrow ~\protect\alpha $-$\protect\gamma $:}

The two roots of%
\begin{equation*}
5-12\theta +8\theta ^{2}=0
\end{equation*}%
are equal to $\frac{3+i}{4}$ and $\frac{3-i}{4}$. Solving (\ref{the}) for $y$
yields $1+i$ and $\frac{1-i}{2}$ for the former, and $1-i$ and $\frac{1+i}{2}
$ for the latter. Selecting one value from each pair (in complex
conjugates), (\ref{exp}) expands to%
\begin{equation*}
\left( y-1-i\right) \left( y-1+i\right) =y^{2}-2y+2
\end{equation*}%
verifying the \textsf{MA(2)} part of the model. The second possible choice
yields%
\begin{equation*}
\left( y-\frac{1-i}{2}\right) \left( y-\frac{1+i}{2}\right) =y^{2}-y+\frac{1%
}{2}
\end{equation*}%
which provides and alternate solution (a `normalized mirror image' of the
previous one).

The three roots of the denominator are all equal to $\frac{5}{4}.$ Solving
for $y$ always yields $2$ and $\frac{1}{2}$. Taking the smaller value of
each pair, we get%
\begin{equation}
\left( y-\frac{1}{2}\right) ^{3}=y^{3}-\frac{3}{2}y^{2}+\frac{3}{4}y-\frac{1%
}{8}  \label{alp3}
\end{equation}
thus recovering the \textsf{AR(3)} part of the model.

We get $V,$ based on (\ref{zer}), by%
\begin{equation*}
V=\frac{\sigma ^{2}\cdot (1-2+2)^{2}}{\omega (0)\cdot \left( 1-\frac{3}{2}+%
\frac{3}{4}-\frac{1}{8}\right) ^{2}}=\frac{176}{2025}
\end{equation*}
(check).

\subparagraph{Spectral$~\longrightarrow ~$correlogram:}

Replacing $\cos \beta $ by (\ref{cos}), $\omega (\beta )$ can be expressed
in a partial-fraction form as%
\begin{equation*}
-\frac{28}{11(y-2)}-\frac{42}{11(y-2)^{2}}-\frac{60}{11(y-2)^{3}}+\frac{31}{%
44\left( y-\frac{1}{2}\right) }+\frac{12}{44\left( y-\frac{1}{2}\right) ^{2}}%
+\frac{15}{176\left( y-\frac{1}{2}\right) ^{3}}
\end{equation*}%
Based on the last three terms (ignoring the first three), we get%
\begin{eqnarray}
\rho _{k} &=&\frac{31}{44}\cdot \left( \frac{1}{2}\right) ^{k-1}+\frac{12}{44%
}\cdot (k-1)\cdot \left( \frac{1}{2}\right) ^{k-2}+\frac{15}{176}\cdot \frac{%
(k-1)(k-2)}{2}\cdot \left( \frac{1}{2}\right) ^{k-3}=  \label{Rk3} \\
&&\left( 1+\frac{3}{44}\cdot k+\frac{15}{44}\cdot k^{2}\right) \left( \frac{1%
}{2}\right) ^{k}  \notag
\end{eqnarray}
which agrees with (\ref{rhox}).

\subparagraph{Correlogram$~\longrightarrow ~$spectral$:$}

Simplifying%
\begin{equation*}
1+2\sum_{k=1}^{\infty }\left( \frac{1}{2}\right) ^{k}\cos (k\cdot \beta )+%
\frac{3}{22}\sum_{k=1}^{\infty }k\cdot \left( \frac{1}{2}\right) ^{k}\cos
(k\cdot \beta )+\frac{15}{22}\sum_{k=1}^{\infty }k^{2}\left( \frac{1}{2}%
\right) ^{k}\cos (k\cdot \beta )
\end{equation*}%
yields (\ref{sec3}).

\subparagraph{Correlogram $\longrightarrow $ $\protect\alpha $-$\protect%
\gamma $ :}

From (\ref{rhox}) it is immediately obvious that there is one triple $%
\lambda $ root of $\frac{1}{2}.$ (\ref{alp3}) then yields the value of the $%
\alpha $ parameters. Knowing these, we can find (\ref{arp3}).

Since $p=3$ and, based on the $k$ range of the general part of (\ref{Rk3}), $%
q-p=-1,$ $q$ must be equal to $2$.

We now need to solve the three nonlinear equations for $\gamma _{1}$ and $%
\gamma _{2}$ which we get by matching the coefficients of $k^{0},$ $k^{1}$
and $k^{2}$ in%
\begin{eqnarray*}
&&\sum_{j,\ell =0}^{2}\left( 1+\frac{15}{22}\cdot (k+j-\ell )+\frac{3}{22}%
\cdot (k+j-\ell )^{2}\right) \cdot \left( \frac{1}{2}\right) ^{j-\ell }\cdot
\gamma _{i}\gamma _{\ell }\overset{}{=} \\
&&\left( 1+\frac{3}{44}\cdot k+\frac{15}{44}\cdot k^{2}\right) \sum_{j,\ell
=0}^{2}\hat{\rho}_{|j-\ell |}\cdot \gamma _{j}\gamma _{\ell }
\end{eqnarray*}%
This is based on (\ref{V1}) and (\ref{Rk3}), where $\hat{\rho}_{k}$ is
defined in (\ref{arp3}), and with the understanding that $\gamma _{0}=1$.

More explicitly, we get a trivial identity when matching the terms without $%
k,$ and (discarding a multiplicative constant)%
\begin{eqnarray*}
22(1+\gamma _{1}^{2}+\gamma _{2}^{2})+42\gamma _{1}(1+\gamma _{2})+27\gamma
_{2} &=&0 \\
22(1+\gamma _{1}^{2}+\gamma _{2}^{2})+30\gamma _{1}(1+\gamma _{2})-9\gamma
_{2} &=&0
\end{eqnarray*}%
(respectively) when matching the $k$ and $k^{2}$ coefficients. The two
solutions are $\gamma _{1}=-2,$ $\gamma _{2}=2$ and $\gamma _{1}=-1,$ $%
\gamma _{2}=\frac{1}{2}$ (check).

\begin{description}
\item[APPENDIX] 
\end{description}

Introducing $y\equiv \exp (i\beta ),$ which implies that 
\begin{equation*}
\cos (k\beta )=\frac{y^{k}+y^{-k}}{2}\ \ \ \ \ \ \text{and \ \ \ \ \ \ \ \ \ 
}d\beta =\frac{dy}{iy}
\end{equation*}%
the RHS of (\ref{fur}) can be converted to%
\begin{equation}
\frac{1}{2\pi i}\doint \omega (y)\left( \frac{y^{k-1}+y^{-k-1}}{2}\right) dy
\label{int}
\end{equation}%
where the integration is over the \emph{unit circle} centered at $0$
(counterclockwise). The theory of \emph{contour integration} tells us that
the result is the sum of \emph{residues} of the integrand at all \emph{%
singularities} inside this circle. Note that the $\omega $ (being,
originally, a function of 
\begin{equation*}
\cos \beta =\frac{y+\frac{1}{y}}{2}
\end{equation*}%
has the $\omega (y)=\omega (\frac{1}{y})$ symmetry.

Expanding $\omega (y),$ which is a \emph{rational function} of $y,$ in terms
of \emph{partial fractions} makes it quite easy to identify the
singularities, and to find their residues. One should realize that each term
of this expansion will have the form of either%
\begin{equation}
c\cdot y^{m}  \label{power}
\end{equation}%
where $c$ is a constant (different for different terms) and $m$ is an
integer (positive, negative, or zero), or of%
\begin{equation}
\frac{c}{(y-\lambda )^{\ell }}  \label{term}
\end{equation}%
where $c$ and $\lambda \neq 0$ are two constants (potentially complex), and $%
\ell $ is a positive integer (usually equal to $1$). Clearly any such term
with complex $c$ and $\lambda $ will have its \emph{complex conjugate}
counterpart. But, more importantly, due to the $y\leftrightarrow \frac{1}{y}$
symmetry, each term of the type $c\cdot y^{m}$ must be `paired' with the
corresponding $c\cdot y^{-m}$ term (this time, with the \emph{same} $c$),
and similarly every term of type (\ref{term}) - regardless whether $c$ and $%
\lambda $ are complex or real - will have its `dual' \emph{set} of terms,
given by the partial-fraction expansion of 
\begin{equation}
\frac{c}{\left( \frac{1}{y}-\lambda \right) ^{\ell }}=\frac{cy^{\ell }}{%
(1-y\lambda )^{\ell }}=\frac{c\left( -\frac{1}{\lambda }\right) ^{\ell
}y^{\ell }}{\left( y-\frac{1}{\lambda }\right) ^{\ell }}  \label{dual}
\end{equation}%
namely 
\begin{equation*}
c\left( -\frac{1}{\lambda }\right) ^{\ell }\sum_{j=0}^{\ell }\binom{\ell }{j}%
\frac{\lambda ^{-j}}{\left( y-\frac{1}{\lambda }\right) ^{j}}
\end{equation*}%
(the details of which become irrelevant shortly). Note that the absolute
value of each $\lambda $ must be different from $1$ (to represent a
stationary model).

One can prove that the contributions of $cy^{m}$ and $cy^{-m}$ to (\ref{int}%
) are identical (for all positive-integer values of $k$); the same is true
for the (\ref{term}) and (\ref{dual}) pair (see below). This means that we
can divide terms of the partial-fraction expansion of $\omega (y)$ into two
parts:

\begin{enumerate}
\item those which have the form of either (\ref{power}) with $m$ \emph{%
positive} or of (\ref{term}) with $|\lambda |<1,$

\item all their \emph{dual} terms - these consist of negative powers of $y,$
and terms or Type (\ref{term}) with $|\lambda |>1,$ all easily identifiable.
\end{enumerate}

We have excluded the $|\lambda |=1$ possibility which cannot happen, as
explained earlier.

All we have to do in the end is to find the sum of residues of terms of Type
1 and multiply the result by $2$ (totally ignoring terms of Type 2, which
would contribute the same amount - that means we never need to utilize the
last two formulas). The constant term is exceptional, since it is its own $%
y\leftrightarrow \frac{1}{y}$ conjugate; it contributes its residue (no
longer multiplied by $2$).

Note that multiplying by $2$ is the same as replacing%
\begin{equation*}
\left( \frac{y^{k-1}+y^{-k-1}}{2}\right)
\end{equation*}%
in (\ref{int}) by $\left( y^{k-1}+y^{-k-1}\right) $ - that's what we will do
from now on.

The residue of each term of Type 1 (and its dual term) follows from 
\begin{equation*}
\frac{1}{2\pi i}\doint y^{m}\left( y^{k-1}+y^{-k-1}\right) dy=\frac{1}{2\pi i%
}\doint y^{-m}\left( y^{k-1}+y^{-k-1}\right) dy=\delta _{m,k}
\end{equation*}%
when $m$ is a non-negative integer. Note that $\delta _{m,k}$ is the
Kronecker's delta (equal to $1$ when $m=k,$ equal to $0$ otherwise).

The result is quite obvious by realizing that the only power of $y$ which
integrates (in the above sense) to $2\pi i$ is $y^{-1}$ - all the other
integer powers yield zero. It is thus the $y^{m-k-1}$ part of the first
integral and the $y^{-m+k-1}$ part of the second integral which result in $%
\delta _{m,k}$ (the $y^{m+k-1}$ and $y^{-m-k-1}$ terms contribute $0$).

Now we want to show that the residue of each term of Type 2 equals to the
residue of its dual term, namely 
\begin{equation*}
\frac{1}{2\pi i}\doint (y-\lambda )^{-\ell }\left( y^{k-1}+y^{-k-1}\right)
dy=\frac{1}{2\pi i}\doint \frac{\left( -\frac{1}{\lambda }\right) ^{\ell
}y^{\ell }}{\left( y-\frac{1}{\lambda }\right) ^{\ell }}\left(
y^{k-1}+y^{-k-1}\right) dy
\end{equation*}%
(assuming that $|\lambda |<1$), and find their values$.$

The LHS integral has two singularities inside the unit circle: one at $0$
and the other at $\lambda $. The first one contributes (another well-known
result of contour integration)%
\begin{eqnarray*}
\left. \frac{1}{k!}\cdot \frac{d^{k}(y-\lambda )^{-\ell }}{dy^{k}}%
\right\vert _{y=0} &=&\frac{(-\ell )(-\ell -1)...(-\ell -k+1)}{k!}\cdot
(-\lambda )^{-\ell -k} \\
&=&(-1)^{\ell }\binom{\ell +k-1}{k}\lambda ^{-\ell -k}
\end{eqnarray*}%
the second one adds%
\begin{eqnarray*}
&&\left. \frac{1}{(\ell -1)!}\cdot \frac{d^{k}y^{k-1}}{dy^{\ell -1}}%
\right\vert _{y=\lambda }+\left. \frac{1}{(\ell -1)!}\cdot \frac{%
d^{k}y^{-k-1}}{dy^{\ell -1}}\right\vert _{y=\lambda }\overset{}{=} \\
&&\frac{(k-1)(k-2)...(k-\ell +1)}{(\ell -1)!}\cdot \lambda ^{k-\ell }+\frac{%
(-k-1)(-k-2)...(-k-\ell +1)}{(\ell -1)!}\cdot \lambda ^{-k-\ell }\overset{}{=%
} \\
&&\binom{k-1}{\ell -1}\lambda ^{k-\ell }+(-1)^{\ell -1}\binom{\ell +k-1}{%
\ell -1}\lambda ^{-\ell -k}
\end{eqnarray*}%
Clearly, the second term cancels the contribution of the $0$ singularity,
with the net result of 
\begin{equation}
\frac{1}{2\pi i}\doint (y-\lambda )^{-\ell }\left( y^{k-1}+y^{-k-1}\right)
dy=\binom{k-1}{\ell -1}\lambda ^{k-\ell }  \label{main}
\end{equation}

The RHS integral has only the $0$ singularity contributing (the $y=\frac{1}{%
\lambda }$ singularity is outside the unit circle), yielding%
\begin{eqnarray*}
\left( -\frac{1}{\lambda }\right) ^{\ell }\cdot \left. \frac{1}{(k-\ell )!}%
\cdot \frac{d^{k}(y-\frac{1}{\lambda })^{-\ell }}{dy^{k-\ell }}\right\vert
_{y=0} &=& \\
\left( -\frac{1}{\lambda }\right) ^{\ell }\cdot \frac{(-\ell )(-\ell
-1)...(-k+1)}{(k-\ell )!}\cdot \left( -\frac{1}{\lambda }\right) ^{-k} &=&%
\binom{k-1}{k-\ell }\lambda ^{k-\ell }
\end{eqnarray*}%
which is the same as (\ref{main}).

We should mention that in the most common case of $\ell =1,$ this answer
boils down to $\lambda ^{k-1}$ (for $\ell =2$ we are getting $(k-1)\lambda
^{k-2},$ for $\ell =3$ this becomes $\frac{(k-1)(k-2)}{2}\lambda ^{k-3},$
etc.).

\end{document}